\documentclass[12pt,reqno,twoside,fleqn]{amsart}%
\usepackage{amssymb}
\usepackage{palatino,hyperref}
\usepackage[mathcal]{euler}
\usepackage[T1]{fontenc}

\makeatletter
\renewcommand{\section}{\@startsection {section}{1}{\z@}%
                                   {-3.5ex \@plus -1ex \@minus -.2ex}%
                                   {.5\linespacing}%
                                   {\normalfont\scshape\centering}}
\makeatother

\numberwithin{equation}{section}

\newtheorem{thm}{Theorem}

\newtheorem{cor}[thm]{Corollary}
  
\theoremstyle{definition}
\newtheorem{definition}{Definition}

\def\beq#1\eeq{\begin{equation}#1\end{equation}}
\makeatletter
\@ifpackageloaded{euler}
{\newcommand{\Wedge}{\wedge}
 \newcommand{\TM}{T\M}
 \newcommand{\TsM}{T^*\M}}
{\newcommand{\Wedge}{\Lambda}
 \newcommand{\TM}{T\mskip-.5\thinmuskip\M}
 \newcommand{\TsM}{T^*\!\M}}
\makeatother
\newcommand{\A}{\mathcal{A}}
\newcommand{\Ak}{\A_\k}
\renewcommand{\k}{\kappa}
\renewcommand{\AA}{\mathbb{A}}
\newcommand{\C}{\mathcal{C}}
\newcommand{\Cci}{\C^\infty_{\mathrm c}}
\newcommand{\M}{\mathcal{M}}

\newcommand{\co}{\mathbb{C}}
\newcommand{\R}{\mathbb{R}}

\newcommand{\VV}{\mathbb{V}}

\newcommand{\inner}{\mathbin{\raise0.3ex\hbox{$\lrcorner$}}}
\newcommand{\Li}{\mathcal{L}}
\newcommand{\D}{\mathcal{D}}
\newcommand{\Hi}{\mathcal{H}}

\DeclareMathOperator{\dom}{dom}

\newcommand{\T}{\mathbb{T}}
\newcommand{\G}{\Gamma}
\DeclareMathOperator{\SU}{SU}
\DeclareMathOperator{\End}{End}
\newcommand{\g}{\mathfrak{g}}

\DeclareMathOperator{\ad}{ad}
\newcommand{\isom}{\mathrel{\widetilde\longrightarrow}}
\newcommand{\Z}{\mathbb{Z}}
\newcommand{\Gc}{\G_{\mathrm c}}
\newcommand{\Omegac}{\Omega_{\mathrm c}}
\DeclareMathOperator{\tr}{tr}
\newcommand{\bbeta}{\hat\beta}
\newcommand{\bgamma}{\hat\gamma}
\newcommand{\Dabrowski}{D¡browski} %D\c{a}browski

\title{Noncommutative Rigidity}
\author{Eli Hawkins}

\begin{document}
\begin{flushright}
\vspace*{-0.5in}
\begin{tabular}{l}
\textsf{\small SISSA 81/2002/FM}\\
\textsf{\small math.QA/0211203}\\
\end{tabular}
\vspace{0.25in}
\end{flushright}
\maketitle
\begin{center}
\vspace{-4ex}
\emph{\small Scuola Internazionale Superiore di Studi Avanzati }\\
\emph{\small Via Beirut 4, I-34014 Trieste, Italy}\\
{\small mrmuon@mac.com}\\
\end{center}

\begin{abstract}
Using very weak criteria for what may constitute a noncommutative geometry,
I show that a pseudo-Riemannian manifold can only be smoothly deformed into noncommutative geometries if certain geometric obstructions vanish.
These obstructions can be expressed as a system of partial differential equations relating the metric and the Poisson structure that describes the noncommutativity. I illustrate this by computing the obstructions for well known examples of noncommutative geometries and quantum groups. These rigid conditions may cast doubt on the idea of noncommutatively deformed space-time.
\end{abstract}

\section{Introduction}
One plausible way to try and construct examples of noncommutative geometry is to start with an ordinary, commutative manifold and deform it. One can try to construct noncommutative algebras that in some sense approximate the algebra of smooth functions on the manifold, and then to construct noncommutative geometries which approximate the geometry of the original manifold. There has been considerable success with the first step. Techniques of geometric quantization can be applied in many cases to construct a sequence of algebras which approximate the algebra of functions in a very strong sense. In a much weaker sense, the formal deformation quantization constructions of Fedosov \cite{fed} and Kontsevich \cite{kon1} give  noncommutative approximations to any manifold. 

Another motive for considering deformations is physical. There are many reasons to suspect that pseudo-Riemannian geometry might not accurately describe the small scale structure of space-time. Noncommutative geometry is a plausible route toward a better description. However, the fact that pseudo-Rie\-man\-nian geometry \emph{is} a sufficient description of space-time for most purposes, suggests that noncommutativity might be treated as a perturbation. 

If so, then this noncommutativity would be described in the leading order by a Poisson structure. Much optimism about this direction was generated  by Kontsevich's remarkable proof \cite{kon1} that there exists a deformation quantization corresponding to \emph{any} Poisson structure. 

Kontsevich's proof was partly inspired by string theory, and this contributed to interest in possible connections between noncommutative geometry and string theory. In \cite{c-d-s} Connes, Douglas, and Schwarz argued that string theory compactified on a torus with a nonzero $3$-form potential is equivalent to string theory compactified on a noncommutative torus. Various authors have argued that a limit of string theory on flat space-time with a constant $2$-form potential ($B$-field) is described by a noncommutative Yang-Mills theory, with the $B$-field providing the Poisson structure (see \cite{s-w,k-s} and refences therein). Cornalba and Schiappo \cite{c-s} argued that the Poisson field is more generally equal to $B+F$, although they went on to advocate using nonassociative ``algebras''. 

Although there exist many examples of noncommutative deformations of the algebra of functions on a manifold, there are very few examples in which this is extended to a deformation of the geometry. The only examples for which the complete axioms of noncommutative geometry are satisfied are the noncommutative torus and a generalization of this constructed by Connes and Landi \cite{c-l} for any compact Riemannian manifold on which $\T^2$ acts by isometries (but see also \cite{c-dv1,c-dv2}). Of course, the axioms of noncommutative geometry are not set in stone, and \Dabrowski\ and Sitarz have constructed an example which only satisfies some of these axioms. On the other hand, Chakraborty and Pal \cite{c-p1} have constructed geometries on $\SU_q(2)$ which satisfy the axioms but do not correspond to the classical geometry of $\SU(2)=S^3$. All this suggests that deforming geometry is not easy and may not be possible generically.

As I explain here, the structures of geometry cannot be deformed in general.
Specifically, integration can only be deformed nicely if the divergence of the Poisson field vanishes, $\pi^{ij}_{\;|j}=0$. The other conditions involve a type of generalized connection called a ``contravariant connection'' (see Sec.~\ref{Contravariant}). The first order differential calculus of 
$1$-forms and the gradient operator $d:\Cci(\M)\to \Omegac^1(\M)$ can only be deformed if there exists a flat, torsion free contravariant connection. A (pseudo)Riemannian metric can only be deformed if 
\[
0 = K_i^{\;jkl}
\]
where $K$ is the curvature of a certain contravariant connection constructed from the Poisson structure and metric. The first two conditions are completely independent of any specific formulation of noncommutative geometry. The derivation of the last condition is motivated by Connes' formulation of noncommutative geometry, but as I explain, it appears to be much more general.

These conditions are necessary, but not sufficient. For one thing, they only depend upon the description of the deformation to leading order; they do not ask whether something may appear in higher order to obstruct the deformation of geometry. This is not even the most complete set of necessary conditions of this kind, but that is for a future paper. 

I anticipate three likely interpretations of this result:
\begin{enumerate}
\item
 Noncommutative deformations have no relevance to the geometry of the universe. 
\item
The geometry of space-time is described by a noncommutative deformation at some level, and my conditions are physical equations of motion of the Poisson field (or whatever fields determine the Poisson structure). 
\item
Noncommutative deformations (in some sense) are relevant, but my assumptions about their properties are too restrictive. 
\end{enumerate}
In the last case, my results can be taken as a guide to what one should not expect or assume when trying to construct noncommutative deformations. 

In Section \ref{Deformations}, I give the necessary background and define what I mean by a deformation. I review how a Poisson structure is derived from a deformation and review the definition of a contravariant connection. In Section \ref{Integration}, I give the obstruction to deforming integration. In Section \ref{1forms}, I give the obstruction to deforming $1$-forms and the gradient operator. In Section \ref{Sec:Metric}, I give the obstruction to deforming a metric. In Sections \ref{Solutions} and \ref{Examples}, I discuss the structure of solutions of these conditions and show how some examples of noncommutative geometry and Poisson manifolds do or do not satisfy these conditions.

\subsection{Notation}
Throughout, $\M$ will denote a locally compact, smooth manifold. $\C^\infty(\M)$ will denote the space of smooth (infinitely differentiable), $\co$-valued functions on $\M$. $\G(\M,V)$ will denote the space of smooth sections of a vector bundle $V$ over $\M$. $\Omega^p(\M) := \G(\M,\Wedge^p \TsM)$ will denote the space of smooth differential $p$-forms. $\Cci$, $\Gc$, and $\Omegac^*$ will denote the spaces of compactly supported smooth functions, sections, and forms.

Lower case Latin characters $f$, $g$, and $h$ will mostly denote smooth functions, except in Section \ref{Sec:Metric} after which $g$ will denote a metric. Lower case Latin characters $a$, $b$, and $c$ will denote elements of a (possibly noncommutative) algebra. Lower case Greek characters $\sigma$, $\rho$, $\alpha$, $\beta$, $\gamma$ will mostly denote $1$-forms. However, $\pi$ will denote the Poisson bivector field, $\omega$ will denote a symplectic $2$-form, and $\epsilon$ will denote a volume form.

Multivectors (vectors, bivectors, \emph{et cetera}) are sections of the exterior powers $\Wedge^*\TM$ of the tangent bundle. I will use the symbol $\inner$ to denote not only the contraction of a vector into a differential form, but also the contraction of a multivector into a differential form. This is such that, for instance, $(X\wedge Y)\inner\epsilon = Y\inner(X\inner\epsilon) = \epsilon(X,Y,\dots)$.

$\Li_X$ denotes the Lie derivative with respect to a vector field $X\in\G(\M,\TM)$. A covariant connection will be denoted as usual by $\nabla$. In index notation this will also be denoted with a vertical stroke as, e.~g., $v^i_{|j}:= \nabla_jv^i$.  A contravariant connection (see Sec.~\ref{Contravariant}) will be denoted as $\D$ to distinguish it from a Dirac-type operator which will be denoted $D$ in Section \ref{Sec:Metric}. 

The symbol $\#$ will denote a certain map, $\#:\TsM\to \TM$, determined by the Poisson structure (Sec.~\ref{Poisson}).
Square brackets will denote the Lie bracket of vector fields or the bracket of a given Lie algebra. Curly brackets will denote the Poisson bracket and its generalizations. $[a,b]_-:=ab-ba$ will denote a commutator. Finally, $[\,\cdot\,,\,\cdot\,]_\pi$ will denote the Koszul bracket of $1$-forms (Sec.~\ref{Poisson})

\section{Deformations}
\label{Deformations}
Let $\M$ be a smooth manifold. I am interested in a smooth, noncommutative deformation of $\M$. What does this mean? 

Suppose that there exists a one-parameter family of algebras $\Ak$  with $\A_0=\Cci(\M)$. We don't need to assume that $\k$ takes a continuous range of values, only that $\Ak$ is defined for all $\k$ in some subset of $\R$ that is dense at $0$. 

Because deformations of this kind are often discussed in the context of the classical limit of a quantum system, the deformation parameter is often called $\hbar$. However, in general, the deformation parameter is not necessarily Planck's constant and I denote it as $\k$ here to maintain this distinction.

Simply having a collection of algebras is not enough. We must tie them together somehow. Thinking of these $\Ak$'s as algebras of smooth functions on noncommutative spaces, the spaces should fit together smoothly into a larger noncommutative space, a sort of noncommutative cylinder. In this way, we have another algebra $\AA$ which is thought of as the algebra of compactly supported, smooth functions on this larger noncommutative space. 

The algebras $\Ak$ form something like a bundle over the set of values of $\k$, although it may not be locally trivial. Let $\Bbbk$ be the algebra of smooth functions of $\k$. 
The algebra $\AA$ should be a $\Bbbk$-algebra; that is, $\AA$ is a $\Bbbk$-module and multiplication in $\AA$ is $\Bbbk$-linear.
We need just a little local triviality. We should be able to (non-canonically) expand any element of $\AA$ as a power series in $\k$ with coefficients in $\Cci(\M)$. Algebraically, this means that $\k^m\AA/\k^{m+1}\AA \cong \Cci(\M)$ for any positive integer $m$.

It is also possible to consider formal deformations in which $\k$ is only a formal parameter. In that case we should take $\Bbbk = \co[[\k]]$, the algebra of formal power series in $\k$. 

In fact, we will only need series expansions to second order in $\k$. At a minimum, we can take for $\Bbbk$ the algebra $\co[\k]/\k^3$ and for $\AA$ the free $\Bbbk$-module $\Cci(\M)\otimes \Bbbk$ with some noncommutative product. 

The trivial case provides some guidance here. Consider the cylindrical space $X\cong \M\times\R$ and the commutative algebras $\AA:=\Cci(X)$ and $\Bbbk=\C^\infty(\R)$; identify $\M$ with $\M\times\{0\}$. There are many ways of identifying $X$ with $\M\times\R$. If we choose one, then this gives a canonical way of expanding any function $f\in\AA=\Cci(X)$ as a power series in $\k$ (the coordinate on $\R$) with coefficients in $\Cci(\M)$. However, if we do not choose such an identification, then such an expansion is not canonical. The ``operator ordering'' ambiguity in a noncommutative deformation is thus equivalent to a coordinate freedom in this commutative case. Nevertheless, any function $f\in\AA=\Cci(X)$ that vanishes along $\M$ can be uniquely written as a multiple of $\k$, any function that vanishes to second order along $\M$ can be written uniquely as a multiple of $\k^2$, and so on.

The following definition is limited to the properties that I will actually need here.
\begin{definition}
\label{def:deformation}
A \emph{smooth deformation} (to second order) of a manifold $\M$ is given by a commutative $\co$-algebra $\Bbbk$, a $\Bbbk$-algebra $\AA$, an element $\k\in\Bbbk$, and an isomorphism $\Cci(\M) \cong \AA/\k\AA$, such that $\AA/\k^3\AA$ is a free $\co[\k]/\k^3$-module.
\end{definition}

The isomorphism $\Cci(\M) \cong \AA/\k\AA$ is a generalization of the requirement that $\A_0 = \Cci(\M)$. The last condition (freeness) expresses the requirement that any element of $\AA$ can be expanded to second order in powers of $\k$ with coefficients in $\Cci(\M)$. Consequently, $\AA/\k^3\AA$ is isomorphic to $\Cci(\M)\otimes \co[\k]/\k^3$ as a $\co[\k]/\k^3$-module. The isomorphism is not canonical. Nevertheless, $\k\AA/\k^2\AA$ is canonically isomorphic to $\Cci(\M)$; that is, the class of any $\k a\in\AA$ modulo $\k^2$  is naturally identified with $\k f$ for some unique $f\in\Cci(\M)$.

We might also require that $\AA$ be a $*$-algebra. That is, that there exists an antilinear involution ${}^*:\A\to\A$ such that $(ab)^*=b^*a^*$. This would ensure that the Poisson field is real. However, this will not affect the form of my constraints at all, so I will not bother discussing it further.

Definition \ref{def:deformation} encompasses both concrete deformations in which the algebras $\Ak$ actually exist, and formal deformations in which $\k$ is only a formal parameter. The concrete case will be important in motivating my definitions.

\subsection{Poisson Structure}
\label{Poisson}
From any smooth deformation, we can construct a Poisson bracket. 

Use the commutator notation, 
\[
[a,b]_- := ab-ba
\mbox,
\]
and use $\bmod\k$ to denote the equivalence class modulo $\k\AA$. For $f,g\in\Cci(\M)$, let $a,b\in\AA$ such that $f=a\bmod\k$ and $g=b\bmod\k$ and define the Poisson bracket $\{f,g\}\in\Cci(\M)$ by
\[
\k \{f,g\} = i [a,b]_- \bmod\k^2
\mbox.\]
The right hand side is a multiple of $\k$ because $[a,b]_- \bmod \k = [f,g]_- =0$.
This defines $\{f,g\}$ uniquely because of the freeness assumption in Definition~\ref{def:deformation}. It only depends on $f$ and $g$, because if we add a multiple of $\k$ to $a$ or $b$, this only changes $[a,b]_-$ by a multiple of $\k^2$. This is explicitly $\co$-bilinear and antisymmetric in $f$ and $g$. 

The properties of the Poisson bracket derive from associativity in $\AA$ to first and second order in $\k$. The identity 
\[
[ab,c]_- = [a,c]_-b + a[b,c]_-
\] 
for $a,b,c\in\AA$, modulo $\k^2$ gives the Leibniz rule (derivation property),
\[
\{fg,h\} = \{f,h\}g + f\{g,h\}
\mbox.
\]
The Jacobi identity for the commutator, modulo $\k^3$, gives the Jacobi identity for the Poisson bracket,
\[
\{f,\{g,h\}\} + \{g,\{h,f\}\} + \{h,\{f,g\}\} =0
\mbox.
\]

The fact that the Poisson bracket is bilinear, antisymmetric, and a derivation on both arguments implies that it is given by a Poisson bivector field $\pi \in \G(\M,\Wedge^2\TM)$. That is, 
\[
\{f,g\} = \pi(df,dg) := \pi^{ij}df_idg_j
\mbox.
\]
The Jacobi identity is equivalent to a condition on $\pi$ which is expressed most succinctly (and cryptically) in terms of the Schouten-Nijenhuis bracket (see e.~g., \cite{vai2}) as
\begin{subequations}
\label{Jacobi.pi}
\beq
0 = [\pi,\pi]
\mbox,
\eeq
or more explicitly as
\beq
0 = \pi^{ia}\pi^{jk}_{\;|a} + \pi^{ja}\pi^{ki}_{\;|a} + \pi^{ka}\pi^{ij}_{\;|a}
\mbox.
\eeq
\end{subequations}
This is a diffeomorphism invariant differential equation; it does not depend on the choice of torsion-free connection.

The Poisson structure determines a vector bundle homomorphism $\#:\TsM\to \TM$, defined in index notation by, $(\#\sigma)^j := \pi^{ij}\sigma_i$. In general, $\#$ is not an isomorphism, but if it is, then its inverse is given by a symplectic structure. That is, there exists a $2$-form $\omega\in\Omega^2(\M)$ which is closed, $d\omega=0,$ and inverse to $\pi$,
\[
\pi^{ai}\omega_{aj} = \delta^i_j
\mbox.\]
In general, the image of $\#$ gives preferred directions of the tangent bundle which determine a \emph{symplectic foliation}. So called because the restriction of the Poisson structure to a leaf is symplectic.

The Koszul bracket $[\,\cdot\,,\,\cdot\,]_\pi$ of $1$-forms is defined by the properties
\beq
\label{pi.bracket1}
[df,dg]_\pi = d\{f,g\}
\eeq
and
\[
[\sigma,f\rho]_\pi = \#\sigma(f) \rho + f[\sigma,\rho]_\pi
\mbox,\]
where $\#\sigma$ acts on $f$ by a directional derivative.
This satisfies the Jacobi identity and is related to the Lie bracket by,
\beq
\label{algebroid}
\#[\sigma,\rho]_\pi = [\#\sigma,\#\rho]
\mbox.\eeq

The cotangent bundle $\TsM$ with the map $\# :\TsM\to \TM$ and the Koszul bracket is an example of a more general structure called a \emph{Lie algebroid} over $\M$. A more elementary example of a Lie algebroid is the tangent bundle itself with the identity map and the Lie bracket.

\subsection{Contravariant Connections}
\label{Contravariant}
I will use the concept of contravariant connection extensively here. This is actually just a case of the natural concept of a connection with respect to a Lie algebroid. Contravariant connections (or ``contravariant derivatives'') were defined by Vaisman \cite{vai2} and analyzed in detail by Fernandes \cite{fer1}; the concept also occurs in \cite{r-v-w} under a different name. 
\begin{definition}
\label{def:contravariant}
Given a vector bundle $V\to\M$, a \emph{contravariant connection} $\D$ is a linear map from $1$-forms to first order differential operators on $V$ such that 
for any $\sigma\in\Omegac^1(\M)$, $f\in\Cci(\M)$, and $v\in\Gc(\M,V)$,
\[
\D_{f\sigma} v = f\,\D_\sigma v
\]
and
\beq
\label{contra.Leibniz}
\D_\sigma(fv) = \#\sigma(f)\, v + f\,\D_\sigma v
\mbox.
\eeq
Here, $\#\sigma(f)$ means the directional derivative of $f$ by the vector field $\#\sigma$.
\end{definition}

This is very similar to the definition of an ordinary (covariant) connection, except that cotangent vectors have taken the place of tangent vectors. This can also be written using index notation such that $\D_\sigma = \sigma_i\D^i$. The contravariant index on $\D^i$ is the reason for the name.

In many ways, contravariant connections behave much like covariant connections. Many standard definitions, identities, and proofs for covariant connections can be translated to contravariant connections simply by exchanging the roles of tangent and cotangent vectors (or covariant and contravariant indices) and replacing Lie brackets with Koszul brackets. A contravariant connection on a vector bundle $V$ is equivalent to one on the dual bundle $V^*$ and determines a contravariant connection on the matrix bundle $\End V = V\otimes V^*$.

We can mimic the standard definition of curvature and define
\beq
\label{curvature}
K(\sigma,\rho)v := \D_\sigma\D_\rho v - \D_\rho\D_\sigma v - \D_{[\sigma,\rho]_\pi} v
\mbox.
\eeq
Note that the Koszul bracket has taken the place of the Lie bracket. A simple calculation shows that $K(\sigma,\rho)$ is a $\Cci(\M)$-linear operator on $\Gc(\M,V)$ and thus all derivatives of $v$ have cancelled. Similarly, $K(\sigma,\rho)$ is $\Cci(\M)$-linear in $\sigma$ and $\rho$. So, $K$ is really an $\End(V)$-valued bivector --- just as the curvature of a covariant connection is an $\End(V)$-valued $2$-form. I will use the term ``flat'' to mean that a contravariant connection has $K=0$.

In the special case that $V=\TsM$, we can mimic the definition of torsion and define,
\beq
\label{torsion}
T(\sigma,\rho) := \D_\sigma\rho - \D_\rho\sigma - [\rho,\sigma]_\pi
\mbox.
\eeq
In the same way, $T$ turns out to be a type $(2,1)$-tensor --- just as the torsion of a covariant connection is a type $(1,2)$-tensor.

The most significant difference between a covariant and a contravariant connection is that $\D_\sigma$ does not necessarily involve any derivatives. If the $1$-form $\sigma$ is such that $\#\sigma =0$, then $\D_\sigma$ is  simply an operator of multiplication by some section of the bundle of matrices $\End V$. When $\D_\sigma$ does involve derivatives, it only takes derivatives parallel to the symplectic foliation determined by $\pi$.

The simplest examples of contravariant connections are derived from covariant connections. If $\nabla$ is a covariant connection, then 
\beq
\label{D.from.Del}
\D_\sigma := \nabla_{\#\sigma}
\eeq 
defines a contravariant connection. However, this clearly has the property that $\#\sigma=0$ implies $\D_\sigma=0$. As this property does not hold in general, not every contravariant connection is of this form.

\section{Integration}
\label{Integration}
If an $n$-dimensional manifold, $\M$, has a volume form $\epsilon\in\Omega^n(\M)$ (e.~g., a Riemannian volume form) then integration of functions defines a linear map, 
\[
\begin{split}
\tau_0 :\,& \Cci(\M) \to \co \\
& f \mapsto \int_\M f\epsilon
\mbox.
\end{split}
\]
This is trivially a trace, $\tau_0(fg)=\tau_0(gf)$, because $\Cci(\M)$ is commutative. When generalizing to noncommutative algebras, it is natural to require integration to be a trace, that is such that $\tau(ab)=\tau(ba)$.  In Connes' formulation of noncommutative geometry, the generalized integration constructed from a spectral triple is automatically a trace.

Given a smooth deformation of $\M$, we can try to smoothly deform integration to a trace. For a concrete deformation, this means that we should have a trace map $\tau_\k:\Ak\to\co$ for each value of $\k$, such that these together give a smooth function of $\k$ for every $a\in\AA$. More generally, and abstractly, we would like a $\Bbbk$-linear trace $\tau:\AA\to\Bbbk$.

Use the notation $\inner$ for the contraction of a vector or multivector into a differential form. In particular $(X\wedge Y)\inner\epsilon = \epsilon(X,Y,\dots)$ and  $\pi\inner\epsilon$ is the differential $(n-2)$-form $(\pi\inner\epsilon)_{ij\cdots} = \frac12 \pi^{ab}\epsilon_{abij\cdots}$.

\begin{thm}
\label{Trace}
Let $\AA$ be a smooth deformation of $\M$. If there exists a $\Bbbk$-linear trace $\tau : \AA \to \Bbbk$ such that
\[
\int_\M f\epsilon = \tau(a) \bmod \k
\]
for $f=a\bmod\k$, then the Poisson field must satisfy,
\begin{subequations}
\label{divergence}
\beq
\label{divergence1}
0 = d(\pi\inner\epsilon)
\mbox,
\eeq
or equivalently, it has $0$ divergence
\beq
\label{divergence2}
0 = \nabla_j\pi^{ij}
\eeq
for any torsion-free connection such that $\nabla\epsilon=0$.
\end{subequations}
\end{thm}
\begin{proof}
The trace property can be equivalently stated as, $0=\tau[a,b]_-$ for all $a,b\in\AA$. Commutativity of $\Cci(\M)$ already implies that $0=\tau[a,b]_-\bmod\k$. The condition for it to vanish to first order in $\k$ is
\begin{align}
0 &= \tau[a,b]_- \bmod \k^2 \nonumber\\
&= -i\k \int_\M \{f,g\}\,\epsilon
\label{Poisson.trace}
\end{align}
for any $a,b\in\AA$ and $f=a\bmod\k$, $g=b\bmod\k$. However, we can rewrite the integrand as 
\[
\{f,g\}\,\epsilon = df\wedge dg \wedge (\pi\inner\epsilon)
\mbox.
\] 
Because $f$ and $g$ are compactly supported, there are  no boundary terms when we rewrite the condition as
\[
0 = \int_\M df\wedge dg \wedge (\pi\inner\epsilon) 
= \int_\M f\, dg \wedge d(\pi\inner\epsilon)
\mbox.\]
Since this must hold for any $f,g\in\Cci(\M)$, this implies eq.~\eqref{divergence1}. 

In index notation, $d(\pi\inner\epsilon)$ is
\[
d(\pi\inner\epsilon)_{ij\cdots}
= \tfrac{n-1}2 \pi^{ab}\!{}_{|[i} \epsilon_{j\cdots ]ab}
= \pi^{ab}_{\:\;|b} \epsilon_{aij\cdots}
\mbox.
\]
So, eq.~\eqref{divergence2} follows.
\end{proof}

Equation \eqref{Poisson.trace} means that integration gives a ``Poisson trace''. In these terms, this result was already given by Weinstein in \cite{wei2}. This condition was applied to quantization in \cite{b-b-g-w,f-s}.

Interestingly, this same conclusion can also be reached from a completely different hypothesis. The orientation axiom for a real spectral triple \cite{con2} requires the existence of a Hochschild homology class which generalizes the volume form of a commutative manifold. If $\epsilon$ can be smoothly deformed into a Hochschild homology class for $\AA$, then \eqref{divergence} must be satisfied. This is indicative of the interplay among Connes' axioms for noncommutative geometry.

\section{Bimodules}
\label{Bimodules}
Let $V$ be a smooth vector bundle over $\M$. The product of a smooth function with a smooth section of $V$ is again a smooth section of $V$. So, as is well known, the space of sections, $\Gc(\M,V)$ is a module of $\Cci(\M)$. It is also a bimodule, if we define left and right multiplication to be the same. Because $\Cci(\M)$ is commutative, this satisfies the bimodule condition of associativity, $f(vg)=(fv)g$ for any $f,g\in\Cci(\M)$ and $v\in\Gc(\M,V)$.

In Sections \ref{1forms} and \ref{Sec:Metric} we will be interested in deforming $\Gc(\M,V)$ into a bimodule of a deformed algebra. In the case of a concrete deformation, for each value of $\k$ we want an $\Ak$-bimodule, and these should fit together nicely into an $\AA$-bimodule, $\VV$, such that $\VV/\k\VV \cong \Gc(\M,V)$. 

\begin{definition}
\label{def:bimodule}
Given a smooth deformation $\AA$ of $\M$, a \emph{bimodule deformation} of a smooth vector bundle $V$ over $\M$ is an $\AA$-bimodule $\VV$ with a $\Cci(\M)$-bimodule isomorphism $\Gc(\M,V)\cong \VV/\k\VV$ such that $\VV/\k^3\VV$ is a free $\co[\k]/\k^3$-module.
\end{definition}

Since $\VV$ is a bimodule, we can define the commutator of $a\in\AA$ and $\hat{v}\in\VV$ in the same way as in an algebra,
\[
[a,\hat{v}]_- := a\hat{v}-\hat{v} a
\mbox.\]
From this, we can construct the (generalized) Poisson bracket of functions and sections, $\{\,\cdot\,,\,\cdot\,\} : \Cci(\M)\times\Gc(\M,V) \to \Gc(\M,V)$.

The following result is present in different terms in \cite{r-v-w}; see also \cite{bur2}. 
\begin{thm}
\label{Flat}
A bimodule deformation of $V$ determines a flat contravariant connection $\D$ on $V$ such that,
\beq
\label{D.from.V}
\D_{df}v = \{ f,v\} 
\eeq
for any $f\in\Cci(\M)$ and $v\in\Gc(\M,V)$.
\end{thm}
\begin{proof}
This bracket satisfies the same properties as the ordinary Poisson bracket (except for antisymmetry, which is just a matter of notation). It is a derivation on both arguments and satisfies the Jacobi identity. The proof is identical to that for the ordinary Poisson bracket.

The Leibniz identity on the first argument is,
\[
\{fg,v\} = \{f,v\}g + f\{g,v\}
\mbox.\]
This implies that $\{f,v\}$ depends locally and linearly on $df$. So, we can define $\D$ by $\D_{df}v := \{f,v\}$.

The Leibniz identity on the second argument is,
\[
\{f,gv\} = \{f,g\} v + g\{f,v\}
\mbox.\]
In terms of $\D$ this is,
\[
\begin{split}
\D_{df}(gv) &= \{f,g\}v + g \, \D_{df}v \\
&= (\# df)g\, v + g\, \D_{df}v
\mbox.
\end{split}
\]
By $\Cci(\M)$-linearity, this property is still true if we replace $df$ with an arbitrary $1$-form $\sigma\in\Omegac^1(\M)$,
\[
\D_\sigma (fv) = \#\sigma(g) \, v + g \,\D_\sigma v
\mbox.\]
This is precisely the condition for $\D$ to be a contravariant connection.

Using equations \eqref{D.from.V}, \eqref{pi.bracket1}, and \eqref{curvature}, the Jacobi identity can be written as
\[
\begin{split}
0 &= \{ f,\{ g,v\}\} - \{ g,\{ f,v\}\} - \{\{f,g\},v\} \\
& = \D_{df}\D_{dg} v - \D_{dg} \D_{df} v - \D_{d\{f,g\}} v \\
& = \D_{df}\D_{dg} v - \D_{dg} \D_{df} v - \D_{[df,dg]_\pi} v \\
& = K(df,dg) v
\end{split}
\]
and therefore $K=0$. That is, $\D$ is a flat contravariant connection.
\end{proof}

Unless $\M$ is a symplectic manifold, this flatness condition is not nearly so restrictive as flatness of an ordinary connection. The equation $0=\D v$ \emph{cannot} be solved locally for a general flat contravariant connection. Nevertheless, this condition is still quite restrictive. 

Using a contravariant connection, one can construct characteristic classes in Poisson cohomology \cite{fer1} from the curvature. These are simply the images of the conventional characteristic classes by the natural map, $\#^* : H^*(\M) \to H^*_\pi(\M)$.  
The existence of a flat contravariant connection thus implies that any (rational) characteristic class of $V$ must be contained in the kernel of $\#^*$, except in degree $0$. In particular, any characteristic class of the restriction of $V$ to a symplectic leaf must be trivial.

\section{One-Forms}
\label{1forms}
Most notions of a noncommutative space involve some generalization of differential forms. In Connes' formulation, a complex of differential forms can be constructed from the algebra (of ``functions'') and the Dirac operator. Other authors (e.~g., \cite{d-k-m1}) simply posit such a complex as a fundamental structure. A common feature of these is that the exterior derivative should be a derivation (satisfying the Leibniz rule). 

In the case of a concrete deformation, suppose that for each value of $\k$ there exists $\Omega^1_\k$ and a derivation $d : \Ak\to\Omega^1_\k$. For $\k=0$ these should be $\Omega^1_0=\Omegac^1(\M)$ and the gradient operator. The Leibniz rule is simply,
\beq
\label{1derivation}
d(ab) = da\,b+a\,db
\mbox.
\eeq
This is only meaningful if $\Omega^1_\k$ is an $\Ak$-bimodule. We should then require the bimodules and derivations to fit together smoothly. This suggests the following definition.
\begin{definition}
\label{def:1forms}
Given a smooth deformation of $\M$, a \emph{deformation of $1$-forms} is a bimodule deformation $\Omega^1$ of the cotangent bundle $\TsM$ along with a derivation
\[
d : \AA \to \Omega^1
\]
that reduces to the gradient operator $d:\Cci(\M)\to\Omegac^1(\M)$ modulo $\k$.
\end{definition}

\begin{thm}
\label{Torsion}
A deformation of $1$-forms determines a contravariant connection on $\TsM$ which is not only flat, but also torsion-free.
\end{thm}
\begin{proof}
By Theorem \ref{Flat}, the bimodule deformation of $\TsM$ determines a flat contravariant connection $\D$ on $\TsM$.

The Leibniz rule \eqref{1derivation} implies that 
\[
d[a,b]_- = [da,b]_- + [a,db]_- = [a,db]_- - [b,da]_-
\mbox,\]
and the Poisson bracket inherits this as,
\[
d\{f,g\} = \{ f,dg\} - \{ g,df\}
\mbox.\]
In terms of the contravariant connection,
\[
\begin{split}
0 &= \D_{df}dg - \D_{dg}df - d\{f,g\} \\
&= \D_{df}dg - \D_{dg}df - [df,dg]_\pi \\
&= T(df,dg)
\mbox.
\end{split}
\]
Therefore $T=0$.
\end{proof}

This is philosophically unsurprising if we consider that torsion of a covariant connection is what happens if we forget that the tangent bundle is the tangent bundle. In this case, the gradient operator establishes that the cotangent bundle really is the cotangent bundle.

In index notation, the property that a contravariant connection is flat and torsion-free can be stated simply as $\D^i\D^j=\D^j\D^i$. 

This is not the only condition of this kind. In fact, the existence of a consistent deformation of $2$-forms places an additional condition on this same contravariant connection. However, the geometric interpretation of this condition is much more complicated and I will postpone discussion to a future paper.

\begin{cor}
\label{Symplectic1}
If $\M$ is a symplectic manifold and there exists a deformation of $1$-forms, then $\M$ admits a flat, torsion-free covariant connection.
\end{cor}
\begin{proof}
Symplectic means that the bundle homomorphism $\#:\TsM\to \TM$ is an isomorphism. In this case, the contravariant connection $\D$ is equivalent to a covariant connection $\nabla$ related by the formula
\beq
\label{intertwine}
\nabla_{\# \sigma} (\#\rho) = \#(\D_\sigma \rho)
\mbox.
\eeq
Note that this differs from eq.~\eqref{D.from.Del}. As noted in \cite{fer1} (Remark 2.3.2) when a covariant and a contravariant connection are intertwined in this way, their curvatures and torsions are intertwined by $\#$. So, the fact that $\D$ is flat and torsion-free implies that $\nabla$ is flat and torsion-free.
\end{proof}

\section{Metric Structure}
\label{Sec:Metric}
The most restrictive (and  concretely defined) notion of noncommutative geometry is given by Connes' axioms for a real spectral triple \cite{con2}. These describe a noncommutative generalization of a Riemannian Spin manifold. 

A real spectral triple involves the following structures: a $*$-algebra $\A$, a Hilbert space $\Hi$ on which $\A$ is faithfully represented, an unbounded self-adjoint operator $D$ on $\Hi$, an antiunitary operator $J$ on $\Hi$, and (for even dimensions) a $\Z_2$-grading on $\Hi$.

For a compact Spin manifold, the algebra is $\A=\C^\infty(\M)$, the Hilbert space is that of square-integrable sections of the spinor bundle, $D$ is the Dirac operator, $J$ is the charge-conjugation operator, and the grading is that into left and right handed spinors.

The axioms for a real spectral triple give the following properties, among others: $\Hi$ is an $\A$-bimodule with the left and right multiplications intertwined by $J$. The common domain, 
\[
\Hi_\infty := \bigcap_{m=1}^\infty \dom D^m
\]
of all powers of $D$ is a projective right (and left) $\A$-module. For any $a\in\A$, the commutator $[D,a]_-$ is a bounded operator on $\Hi$ (or $\Hi_\infty$) and commutes with right-multiplication by any $b\in\A$.

In the case of a compact Spin manifold, $\Hi_\infty$ is the space of smooth sections of the spinor bundle and $[D,f]_- = i\gamma^i df_i$, in terms of the Dirac matrix-vector.

For a spectral triple, there is a construction of noncommutative differential forms which is quite simple in degree $1$. The $\A$-bimodule $\Omega^1_D\subset\Li(\Hi)$ is generated by all commutators $[D,a]_-$ for $a\in\A$ and the differential map $d:\A\to\Omega^1_D$ is given by $da = -i[D,a]_-$.

The properties that we will need here are more general than Spin manifolds and real spectral triples. We do not need the operator $J$, only the bimodule structure. We do not need the Hilbert space $\Hi$, only the bimodule $\Hi_\infty$. Indeed, $\Hi$ could even be a Krein space, as in the approach to noncommutative space-time advocated by Strohmaier \cite{stro}. We do not need spinors, only some Spin bundle; that is, a bundle carrying a representation of the bundle of Clifford algebras. The Spin bundle $\Wedge^*\TsM$, which exists for any Riemannian or pseudo-Riemannian manifold, is sufficient.

\begin{definition}
A \emph{deformed spectral triple} $(\AA,\VV,D)$ for a (pseudo)Riemannian manifold $\M$ consists of a smooth deformation $\AA$ of $\M$, a bimodule deformation $\VV$ of some Spin bundle $V$, and a $\Bbbk$-linear operator $D:\VV\to\VV$ such that
\begin{enumerate}
\item
$D \bmod \k$ is the Dirac-type operator on sections of $V$.
\item
For any $a\in\AA$, $[D,a]_-$ commutes with the right multiplication of $\AA$ on $\VV$.
\item 
The construction of noncommutative differential forms from $\AA$ and $D$ gives a deformation of $1$-forms on $\M$.
\end{enumerate}
\end{definition}

Let $\Omega^1:=\Omega^1_D$ be the space of noncommutative $1$-forms constructed with $\AA$ and $D$.
Modulo $\k$, this reduces to the classical space of $1$-forms, $\Omegac^1(\M)$. The isomorphism $c : \Omegac^1(\M) \isom \Omega^1/\k\Omega^1$ is the Clifford representation on $V$; if $V$ is the spinor bundle, this is given by the Dirac matrix-vector as $c(\sigma) = \gamma^j\sigma_j$. The only part of the definition of a deformation of $1$-forms which is not automatically satisfied is the condition that $\Omega^1/\k^3\Omega^1$ be a free $\co[\k]/\k^3$-module. This is effectively a condition on the behavior of $D$. This condition is automatically satisfied for an ``isospectral'' deformation in which $\VV$ can be represented as a free $\Bbbk$-module such that $D$ is independent of $\k$. 

\begin{thm}
\label{Metric}
\label{def:triple}
A deformed spectral triple determines a flat, torsion-free contravariant connection on $\TsM$ which is compatible with the metric, 
\beq
\label{metric}
0 = \D^i g^{jk}
\mbox.
\eeq
\end{thm}
\begin{proof}
Since by hypothesis $\VV$ and $\Omega^1$ are bimodule deformations of $V$ and $\TsM$, we have flat contravariant connections on $V$ and $\TsM$ by Theorem \ref{Flat}. By Theorem \ref{Torsion}, the contravariant connection on $\TsM$ is also torsion-free.

For, $a\in\AA$, $\hat\sigma\in\Omega^1$, and $\hat{v}\in\VV$, associativity gives the elementary identity,
\beq
\label{spinor.com}
[a,\hat\sigma\,\hat{v}]_- = [a,\hat\sigma]_-\hat{v} + \hat\sigma [a,\hat{v}]_-
\mbox.
\eeq
This is automatically true to $0$'th order in $\k$. At first order, \eqref{spinor.com} gives an identity for generalized Poisson brackets; let $f\in\Cci(\M)$, $\sigma\in\Omegac^1(\M)$, and $v\in\Gc(\M,V)$, 
\[
\{f,c(\sigma)v\} = c(\{f,\sigma\})v + c(\sigma) \{f,v\}
\mbox.
\]
In terms of the contravariant connections for $\TsM$ and $V$,
\[
\D_{df}[c(\sigma)v] = c(\D_{df}\sigma)v + c(\sigma)\D_{df}v
\mbox,\]  
and so $c$ is ``parallel'' with respect to the contravariant connections.

The Clifford identity,
\[
g^{ij}\sigma_i\rho_j = \tfrac12\left[c(\sigma)c(\rho) + c(\rho)c(\sigma)\right]
\]
then shows that the metric pairing is also parallel, which is equivalent to \eqref{metric}.
\end{proof}

Other arguments are possible for this same metric-compatibility condition \eqref{metric}. For instance, given a deformation of $1$-forms, we might try to define a deformed noncommutative metric as an $\AA$-bimodule homomorphism $\langle\,\cdot\,,\,\cdot\,\rangle : \Omega^1\otimes_\AA\Omega^1 \to \AA$ which reduces modulo $\k$ to the classical metric (in contravariant form). By bimodule linearity, for $a\in\AA$ and $\bbeta,\bgamma\in\Omega^1$,
\[
[a,\langle\bbeta,\bgamma\rangle]_- = \langle[a,\bbeta]_-,\bgamma\rangle + \langle\bbeta,[a,\bgamma]_-\rangle
\mbox.
\]
Modulo $\k$, this gives
\[
(\#df)\langle\beta,\gamma\rangle = \langle\D_{df}\beta,\gamma\rangle + \langle\beta,\D_{df}\gamma\rangle
\mbox,
\]
which implies the metric-compatibility property \eqref{metric}.  

Alternatively, suppose that $\mathbb{S}^2\subset\Omega^1\otimes_\AA\Omega^1$ is a bimodule deformation of the symmetric $2$-tensor bundle $S^2\TsM$. For instance, this might be  the kernel of a deformed exterior product map $\wedge : \Omega^1\otimes_{\AA}\Omega^1 \to \Omega^2$. We might  define a deformed noncommutative metric to be a bimodule homomorphism $\mathbb{S}^2 \to \AA$ which reduces to the classical metric modulo $\k$. The deformation determines a contravariant connection on $S^2\TsM$ which coincides with that induced by the connection for $\TsM$. Again, this implies metric compatibility.

\begin{cor}
\label{Symplectic2}
If a Riemannian manifold admits a deformed spectral triple such that the Poisson structure is symplectic, then it admits a flat Riemannian metric defined by 
\beq
\label{g'}
g'_{jk} := \omega_{ja}\omega_{kb}g^{ab}
\mbox.
\eeq
\end{cor}
\begin{proof}
By Corollary \ref{Symplectic1}, the contravariant connection $\D$ for $\TsM$ is equivalent to a flat, torsion-free covariant connection $\nabla$. This is related by eq.~\eqref{intertwine}. In terms of index notation, we can translate between expressions in terms of $\D$ and $\nabla$ by raising and lowering all indices with the Poisson field $\pi$ and the symplectic form, respectively.

Because $\D$ is compatible with the metric $g$, we have
\[
\begin{split}
0 &= \omega_{ia}\omega_{jb}\omega_{kc} \D^ag^{bc}\\
&= \nabla_i (\omega_{jb}\omega_{kc}g^{bc}) \\
&= \nabla_i g'_{jk}
\mbox.
\end{split}
\]
Therefore, $\nabla$ is the Levi-Civita connection for $g'$ and its curvature, $0$, is the Riemannian curvature of $g'$.
\end{proof}
Note that if $g$ is Riemannian, then so is $g'$. In general, $g'$ has the same signature as $g$.

\emph{A priori}, the condition presented by Thm.~\ref{Metric} is existential. One might despair of proving the nonexistence of a suitable connection. However, the situation is actually much better. There is a strong formal analogy between formulas in Riemannian geometry and in the present setting; the roles of covariant and contravariant indices are simply interchanged. In particular, there is an analogue of the Levi-Civita connection.

\begin{thm}
\label{Levi-Civita}
If a manifold $\M$ has both a metric $g$ and a Poisson field $\pi$, then there exists a unique torsion-free, metric-compatible, contravariant connection $\D$  given by,
\begin{multline}
\label{metric.connection}
\langle\D_\alpha\beta,\gamma\rangle
= \tfrac12 \bigl[(\#\alpha) \langle\beta,\gamma\rangle - (\#\gamma) \langle\alpha,\beta\rangle +  (\#\beta) \langle\gamma,\alpha\rangle \\
  + \langle[\gamma,\alpha]_\pi,\beta\rangle - \langle[\beta,\gamma]_\pi,\alpha\rangle + \langle[\alpha,\beta]_\pi,\gamma\rangle\bigr]
\mbox,
\end{multline}
for $\alpha,\beta,\gamma\in \Omega^1(\M)$. Here, $\langle\alpha,\beta\rangle:= g^{ij}\alpha_i\beta_j$ is the metric pairing and the vectors $\#\alpha$, \emph{et cetera} act by directional derivatives.
\end{thm}
\begin{definition}
This $\D$ is the \emph{metric contravariant connection}.
\end{definition}
\begin{proof}
The construction is precisely analogous to that of the Levi-Civita connection. 

Suppose that such a connection exists. The metric compatibility condition implies that
\[
(\#\alpha) \langle\beta,\gamma\rangle = \langle\D_\alpha\beta,\gamma\rangle + \langle\beta,\D_\alpha\gamma\rangle
\mbox,\]
and the torsion-free condition is
\[
\D_\alpha\gamma = \D_\gamma\alpha - [\gamma,\alpha]_\pi 
\mbox.\]
Putting these facts together gives
\[
\langle\D_\alpha\beta,\gamma\rangle
= (\#\alpha) \langle\beta,\gamma\rangle + \langle[\gamma,\alpha]_\pi,\beta\rangle - \langle\D_\gamma\alpha,\beta\rangle 
\mbox.
\]
By iterating this three times, we can solve for $\langle\D_\alpha\beta,\gamma\rangle$ and get eq.~\eqref{metric.connection}.

This construction proves uniqueness. To prove that this defines a connection, one can multiply the arguments by functions to show that all derivatives of $\alpha$ and $\gamma$ cancel and verify that this satisfies the correct Leibniz rule in the argument $\beta$. 
\end{proof}

The formula \eqref{metric.connection} is the natural analogue of a formula for the Levi-Civita connection. In fact \eqref{metric.connection} has already appeared, although without derivation, in \cite{bou}.

\begin{cor}
\label{Metric2}
If a (pseudo)Riemannian manifold $\M$ admits a deformed spectral triple, then the metric contravariant connection must be flat.
\end{cor}
\begin{proof}
This is immediate from Theorems \ref{Metric} and \ref{Levi-Civita}.
\end{proof}

This means that the necessary condition for compatibility of Riemannian and Poisson structures may be expressed as a differential equation, albeit a very cumbersome one when written out explicitly.

We can also construct the metric contravariant connection in terms of the Levi-Civita (covariant) connection, $\nabla$. The contravariant connection $\nabla_{\#}$ defined by eq.~\eqref{D.from.Del} is metric-compatible, but  has torsion. Using equations \eqref{pi.bracket1} and \eqref{torsion}, the torsion is
\[
T^{\nabla_\#}(df,dh) = \nabla_{\#df}dh - \nabla_{\#dh}df - d\{f,h\}
\mbox.
\]
Using index notation (with a vertical stroke for $\nabla$) this is,
\[
T^{ij}_{\;k} f_{|i}h_{|j} = \pi^{ij}f_{|i} h_{|jk}+ \pi^{ij}f_{|ik}h_{|j} - (\pi^{ij}f_{|i}h_{|j})_{|k}
= -\pi^{ij}_{\; |k} f_{|i}h_{|j}
\mbox.
\]
Thus $T^{ij}_{\;k} = -\pi^{ij}_{\;|k}$ or equivalently,
\[
T^{\nabla_\#} = -\nabla\pi
\mbox.\]
The metric contravariant connection can be constructed by correcting this; its action on a $1$-form is
\[
\D^i\sigma_k = \pi^{ij}\nabla_j\sigma_k + A^{ij}_{\,k}\sigma_j
\mbox,\]
where 
\beq
\label{A.form}
A^{ij}_{\,k} := \tfrac12 \left(\pi^{ij}_{\;|k} - \pi^{j\;\: |i}_{\;k} - \pi^{i\;\: |j}_{\;k}\right)
\eeq
and indices are raised and lowered with the metric. 

Using this, we can explicitly construct the curvature of $\D$. It is
\beq
\label{explicit.K}
K^{ijk}{}_l = \pi^{ia}\pi^{jb}R^k{}_{l ab} - \pi^{ia}A^{jk}_{\; l \; | a} + \pi^{ja} A^{ik}_{\; l\;|a} - A^{ia}_{\;l} A^{jk}_{\;a} + A^{ja}_{\;l} A^{ik}_{\;a} + \pi^{ij}_{\;|a} A^{ak}_{\;\:l}
\mbox.\eeq
Although it is not at all apparent from this expression, this tensor actually has precisely the symmetries of the Riemann tensor. Setting this to $0$ gives a differential equation relating $\pi$ and the metric.

\section{Some Solutions}
\label{Solutions}
If a commutative manifold is deformed into noncommutative geometry, it ought to have both a metric structure and integration. If a manifold can be deformed with respect to a given Poisson structure, then we should expect the conclusions of both Theorem \ref{Trace} and Corollary \ref{Metric2} to hold. That is, it should satisfy both the conditions of $0$ divergence and flatness of the metric contravariant connection. There may indeed be other necessary conditions, but just considering the conditions at hand is illuminating.

There is one case in which a general solution to the second condition can be written down explicitly. Suppose that the manifold $\M$ is $\R^{2n}$ and the Poisson field is symplectic. Because $\R^{2n}$ is contractible, the symplectic form must be exact and can be written in terms of a potential $1$-form as $\omega=d\theta$. According to Cor.~\ref{Symplectic2}, a flat metric $g'$ can be constructed from $\omega$ and the original metric $g$. We can spend our coordinate freedom and  fix $g'$ to be some constant metric on $\R^{2n}$. The metric $g$  and symplectic form are now both determined by $\theta$ as 
\beq
\label{g.from.omega}
g_{ij} = \omega_{ia}\omega_{jb}{g'}^{ab}
\eeq
 and $\omega=d\theta$.

Let $\epsilon$ and $\epsilon'$ be the volume $2n$-forms of the metrics $g$ and $g'$, respectively. Since all volume forms are proportional, the symplectic volume form can be written as 
\beq
\label{symplectic.volume}
\frac{\omega^n}{n!} = h \epsilon'
\eeq
 in terms of some smooth function $h\in\C^\infty(\R^{2n})$. In coordinates the volume forms are given by the determinants of $g$, $g'$, and $\omega$. Equation \eqref{g.from.omega} gives 
\[
\det g = (\det \omega)^2(\det g')^{-1}
\]
and eq.~\eqref{symplectic.volume} gives $\det \omega = h^2\det g'$. Thus
\[
(\det\omega)^2 = \det g \det g' = h^4 (\det g')^2
\]
and $\epsilon = h^2\epsilon' = h \frac{\omega^n}{n!}$. We can therefore write, $\pi\inner\epsilon = h \frac{\omega^{n-1}}{n!}$ and the divergence condition is
\[
0 = d(\pi\inner\epsilon) = dh \wedge \frac{\omega^{n-1}}{n!}
\mbox;
\]
That is, that $h$ is constant.

So, the general solution of the conditions for symplectic $\R^{2n}$ is given by a constant metric $g'$ on $\R^{2n}$ and a $1$-form $\theta\in\Omega^1(\R^{2n})$ such that $\det d\theta$ is constant. There is a gauge freedom to add any gradient to $\theta$, thus there are $2n-2$ local degrees of freedom to this solution. 

This is in contrast to an arbitrary Riemannian metric on $\R^{2n}$. That has $n(2n+1)$ components, but there are $2n$ degrees of diffeomorphism freedom per point. Therefore an arbitrary Riemannian metric has $n(2n-1)$ degrees of freedom per point.

Of course, an arbitrary symplectic manifold resembles this locally, so the counting of degrees of freedom is not limited to $\R^{2n}$. In two dimensions this is particularly striking. There are no local degrees of freedom. In this way, only a flat $2$-torus can be deformed into a noncommutative geometry.

From this it is clear that only a very restricted class of Riemannian geometries is compatible with any symplectic structure. The requirement of compatibility with a deformation leaves geometry quite rigid.
\medskip

The condition of contravariant flatness given by Corollary \ref{Metric2} can be considered as a second order nonlinear partial differential equation in $g$ and $\pi$. It is homogeneously quadratic in $\pi$. Notably, it only involves derivatives in directions in which $\pi$ is nonvanishing; that is, along the symplectic foliation determined by $\pi$. This gives a sort of causal structure. We should expect solutions to propagate only along the symplectic leaves.

To get a greater sense of the structure of solutions of these conditions, we can consider perturbations of a given solution. However, perturbations of $\pi$ are in general badly behaved. If the rank of $\pi$ (in some neighborhood) is less than the dimension of $\M$, then an arbitrarily small perturbation can increase the rank. This would increase the dimension of the symplectic leaves and thus drastically change the causal structure of the contravariant flatness condition. 

It is thus much easier to consider perturbations of the metric relative to a fixed Poisson structure. The curvature of the metric contravariant connection is in some ways analogous to the curvature of the Levi-Civita connection, so consider the analogous problem: perturbing a flat metric to other flat metrics. If we are given a flat metric on some manifold, then locally any other flat metric is equivalent via a diffeomorphism. So, locally, any perturbation of the flat metric to another flat metric is given by an infinitesimal diffeomorphism; that is, the perturbation is the Lie derivative with respect to some vector field $\xi$,
\beq
\label{Lie.g}
\delta g_{ij} = \Li_\xi g_{ij} = \nabla_i\xi_j + \nabla_j\xi_i
\mbox,
\eeq
where the index on $\xi_i$ is of course lowered with the metric.

The linearized contravariant flatness condition is formally very similar. Let $h^{ij} = \delta g^{ij}$ be a perturbation of the contravariant form of the metric, and suppose that the metric contravariant connection of $g$ is flat. The condition that the perturbation preserve flatness is,
\beq
\label{linearized}
0 = \D^j\D^l h^{ik} + \D^i\D^k h^{jl} - \D^k\D^l h^{ij} - \D^i\D^j h^{kl}
\mbox.
\eeq
Note that the order of the $\D$'s is irrelevant because $\D$ is flat and torsion-free. Using the Koszul bracket in place of the Lie bracket, we can define a sort of Lie derivative with respect to a $1$-form. This can be expressed in terms of $\D$, and the analogue of \eqref{Lie.g} is, for some $\alpha\in\Omega^1(\M)$,
\beq
\label{pi.Lie.g}
h^{ij} = \D^i\alpha^j + \D^j\alpha^i
\mbox.
\eeq
It is easy to check that this is a solution of eq.~\eqref{linearized}. In fact, the computation is formally identical to checking that \eqref{Lie.g} is a solution of the linearized Riemann flatness condition.

The formula \eqref{pi.Lie.g} gives enough solutions locally in the symplectic case, but not generically. For example, if $\pi=0$, then any metric satisfies the contravariant flatness condition; so, any perturbation will satisfy eq.~\eqref{linearized}. Yet \eqref{pi.Lie.g} is identically $0$. Clearly, there are other solutions in general.

If the metric is perturbed by $h^{ij}$, then the change in the volume form is $\delta\epsilon = -\frac12 h^k{}_k \epsilon$. The linearization of eq.~\eqref{divergence1} is thus,
\[
0 = d(h^k{}_k \pi\inner\epsilon) = d(h^k{}_k)\wedge (\pi\inner\epsilon)
\]
or simply, $0= \pi^{ij}\nabla_j h^k{}_k$. That is, $h^k{}_k$ must be constant along the symplectic leaves.

Inserting the expression \eqref{pi.Lie.g} into this equation gives the condition,
\beq
\label{alpha}
\begin{split}
0 &= \pi^{ij}\nabla_j (\D^k\alpha_k) \\
&= \#d(\D^k\alpha_k)
\mbox.
\end{split}
\eeq
That is, $\D^k\alpha_k$ should be constant along symplectic leaves.
So, a class of perturbations preserving the conditions for deformation is given by \eqref{pi.Lie.g} for $\alpha$ satisfying eq.~\eqref{alpha}.

\section{Examples and Counterexamples}
\label{Examples}
\subsection{The Noncommutative Torus}
The oldest and best known example of a noncommutative geometry obtained by deforming a commutative manifold is the noncommutative $2$-torus. The real spectral triple which describes the geometry was given in \cite{con2}. The classical Riemannian manifold is a flat $2$-torus. The classical geometry and the deformed algebra are both invariant under the action of the $2$-torus group $\T^2$, therefore the Poisson field must be invariant. The Poisson structure is thus symplectic with symplectic form a multiple of the volume form. The Levi-Civita connection is just the trivial connection, so
\[
0 = \nabla_k\pi^{ij}
\]
and in particular, $\nabla_j\pi^{ij}=0$ which is precisely what we should expect from Thm.~\ref{Trace}. The metric contravariant connection must be $\T^2$ invariant, so the simplest possibility is the connection $\nabla_\#$ derived from the Levi-Civita connection. A simple computation shows that this is flat, torsion-free, and metric-compatible; thus the conclusion of Thm.~\ref{Metric} also holds.

\subsection{Isospectral Deformations}
The noncommutative torus was generalized considerably by Connes and Landi \cite{c-l}. The construction applies to a compact Riemannian Spin manifold $\M$ on which a torus group $\T^m$ ($m\geq2$) acts isometrically. We can assume without loss of generality that $\T^m$ acts faithfully.

Given any constant (i.~e., $\T^m$-invariant) Poisson field on $\T^m$, we can deform the algebra of smooth functions to a noncommutative algebra $\C^\infty(\T^m_\k)$ which is a direct generalization of the noncommutative $2$-torus. The algebra $\C^\infty(\M)$ is then twisted by the noncommutative $\T^m$. The deformed algebra is the $\T^m$-invariant subalgebra of $\C^\infty(\M)\mathbin{\widehat\otimes} \C^\infty(\T^m_\k)$.

The construction of the spectral triple is quite simple. The deformed algebra is represented on the same Hilbert space of square-integrable sections of the classical spinor bundle. The classical Dirac operator is used with the deformed algebras as well. Since the operators are identified, their spectra are the same and this gives the name ``isospectral'' to this construction.

This family of spectral triples is very well behaved. They satisfy the axioms \cite{con2} for real spectral triples. In particular, a trace on each algebra can be constructed using the Dixmier trace and the Dirac operator. We should certainly expect that these examples will fit my definitions and satisfy the conditions I have formulated. 

Let $\{X_a\}$ be a basis of infinitesimal generators of the $\T^m$-action on $\M$ and let $\Pi^{AB}$ be the components of the (constant) Poisson field on $\T^m$. The Poisson field on $\M$ is simply $\pi = \frac12 \Pi^{AB} X_A\wedge X_B$.

The vectors $X_A$ are Killing vectors, so in particular
\[
0 = \Li_{X_A}\epsilon = d(X_A\inner\epsilon)
\mbox.\]
Because $\T^m$ is abelian, $[X_A,X_B]=0$, and thus
\[
\begin{split}
0 &= [X_A,X_B]\inner\epsilon \\
&= \Li_{X_A}(X_B\inner\epsilon) - X_B\inner(\Li_{X_A}\epsilon) \\
&= X_A\inner d(X_B\inner\epsilon) + d\left[X_A\inner(X_B\inner\epsilon)\right] \\
&= -d\left[(X_A\wedge X_B)\inner\epsilon\right]
\mbox.
\end{split}
\]
This implies that $d(\pi\inner\epsilon)=0$, so this indeed satisfies the conclusion of Thm.~\ref{Trace}.

The map $\#:\TsM\to \TM$ is given by, $\#\sigma = \Pi^{AB}(X_A\inner\sigma)X_B$. Define a contravariant connection by $\D_\sigma = \Pi^{AB} (X_A\inner\sigma) \Li_{X_B}$. Because the $X_A$'s are Killing vectors, we have immediately metric-compatibility, $\D g=0$.

It is sufficient to compute the torsion with exact $1$-forms. So for $f,h\in\C^\infty(\M)$ we find the derivative
\[
\D_{df} dh = \Pi^{AB} X_A(f)\, d(X_Bh)
\mbox.
\]
This gives
\[
\begin{split}
T(df,dh) &= \D_{df}dh - \D_{dh}df - d\{f,h\} \\
&= \Pi^{AB}\left[X_Af\,d(X_Bh) + d(X_Af)\,X_Bh - d(X_Bf\,X_Bh)\right] \\
&= 0
\mbox,
\end{split}
\]
and hence $T=0$. Therefore, by the uniqueness result of Thm.~\ref{Levi-Civita}, $\D$ is the metric contravariant connection.

If $\sigma\in\Omega^1(\M)$ is $\T^m$-invariant, then $\D\sigma=0$, so clearly $K\sigma =0$. This shows that $K=0$ at any point that is not fixed by any proper subgroup of $\T^m$, but such points are dense and therefore $K=0$ identically. So, with this contravariant connection, $\M$ satisfies the conclusion of Thm.~\ref{Metric}.

\subsection{The Fuzzy Sphere}
The ``fuzzy sphere'' \cite{mad1} is another popular noncommutative space. It is a concrete, $\SU(2)$-equivariant deformation of the sphere $S^2$. Algebras $\Ak$ exist for $\k=0$ and $\k = N^{-1}$ ($N$ a positive integer).
 The algebra $\A_{1/N}$ is simply the algebra of $N$ by $N$ matrices. The equivariance provides a way of simulating geometric constructions on the fuzzy sphere, but a spectral triple that describes its geometry has not been found.

Because the deformation is equivariant, the Poisson field $\pi$ must be $\SU(2)$-in\-vari\-ant. The only possibility is symplectic, with symplectic form some multiple of the volume form. The contraction $\pi\inner\epsilon$ is a constant function on $S^2$, so $d(\pi\inner\epsilon)=0$ and the conclusion of Thm.~\ref{Trace} is satisfied. Indeed, integration on $S^2$ can be smoothly deformed to traces. On the algebra $\A_{1/N}$, this is simply $\tau_{1/N} = \frac{4\pi}N \tr$.

By Cor.~\ref{Symplectic2}, there can be no deformed spectral triple for the fuzzy sphere, because $S^2$ does not admit a flat metric. By Cor.~\ref{Symplectic1}, it is not even possible to deform $1$-forms to the fuzzy sphere, because $T^*S^2$ does not admit any flat connection.

The fuzzy sphere generalizes nicely to other coadjoint orbits. A coadjoint orbit of a compact, semisimple Lie group $G$ is deformed to a sequence of $G$-equivariant matrix algebras. The Poisson structures are symplectic and $G$-invariant. None of these admit deformed spectral triples, because none of these coadjoint orbits admit a flat metric.
 
\subsection{Two Dimensions}
Suppose that $\M$ is a $2$-dimensional, connected Poisson manifold with a volume form $\epsilon\in\Omega^2(\M)$, and let $\epsilon^{ij}$ be the inverse of the volume form. Because the bundle of bivectors is rank $1$, we must have $\pi^{ij} = h \epsilon^{ij}$ for some smooth function $h\in\C^\infty(\M)$. 

This ansatz gives that $\pi\inner\epsilon = h$. So, the condition from Thm.~\ref{Trace} that $0=d(\pi\inner\epsilon) = dh$ implies that $h$ is constant. In other words, either $\pi=0$ or $\M$ is symplectic with symplectic form a constant multiple of $\epsilon$. 

If $\M$ is also Riemannian and admits a deformed spectral triple, Corollary \ref{Metric2} then shows that if $\pi\neq0$ then $\M$ admits a flat metric. That is, it can only be a torus $\T^2$, a cylinder $S^1\times\R$, or the plane $\R^2$.

Suppose instead that $\M$ admits a deformed spectral triple, but disregard the hypothesis of Thm.~\ref{Trace} (that integration is deformed to a trace). The open submanifold 
\[
\Sigma:= \{x\in\M\mid \pi(x)\neq0\}
\] 
is symplectic with symplectic form $\omega = h^{-1} \epsilon$. If $\M$ admits a deformed spectral triple, then by Cor.~\ref{Symplectic2} there exists a flat metric on $\Sigma$ given by,
\[
g'_{ij} = \omega_{ia}\omega_{jb} g^{ab} = h^{-2} g_{ij}
\mbox.
\]
In other words, $g_{ij}$ must be flattened by the conformal factor $h^{-1}$. Of course, any $2$-dimensional manifold is locally conformally flat, and so such an $h$ always exists locally.

Consider the unit sphere $S^2$. This is not globally conformally flat, but if we remove a point, then $\Sigma = S^2\smallsetminus *$ is. We can write the metric in an explicitly conformally flat form using a complex coordinate $\zeta$,
\beq
\label{S2.metric}
ds^2 = 4 (\zeta\bar\zeta + 1)^{-2}d\zeta\,d\bar\zeta .
\eeq
This shows that the condition will be satisfied by $h = 2c (\zeta\bar\zeta + 1)^{-1}$ for $c$ constant. In terms of the standard embedding coordinates, $h = c (1+z)$, so this is a smooth function on $S^2$ and vanishes at the deleted point.

The volume form for the metric \eqref{S2.metric} is $\epsilon = 2i (\zeta\bar\zeta + 1)^{-2}d\zeta\wedge d\bar\zeta$, and so the symplectic form will be $\omega = h^{-1}\epsilon = -ic^{-1} (\zeta\bar\zeta + 1)^{-1} d\zeta\wedge d\bar\zeta$. The Poisson bracket is thus given by,
\beq
\label{zeta.Poisson}
\{\zeta,\bar\zeta\} = -ci (\zeta\bar\zeta + 1) .
\eeq

The sphere minus two points is also globally conformally flat. However, the conformal factor is such that $h$ would not be differentiable on $S^2$. This would mean that the Poisson bracket of two smooth functions would not always be differentiable. This does not correspond to a smooth deformation.

The Podle\'{s} standard sphere is a (noncommutative) homogeneous space of the quantum group $\SU_q(2)$. For $q=1$ it is the commutative sphere, $S^2$. \Dabrowski\ and Sitarz have constructed an unusual spectral triple for the Podle\'{s} sphere which is $\SU_q(2)$-equivariant and reduces to the spectral triple for a sphere of radius $1$ when $q=1$; it satisfies some --- but not all --- of Connes' axioms for a real spectral triple. 

The Podle\'{s} standard sphere algebra $\A(S^2_q)$ (for $q\in\R$) is the $*$-algebra generated by $A$ and $B$ such that,
\[
\begin{split}
A&=A^*, & 
AB &= q^2 BA, &
BB^* &= q^{-2} A(1-A), &
B^*B &= A(1-q^2A) 
\mbox.
\end{split}
\]
This is a smooth deformation of $S^2$ with $q=1+\k$. When $q=1$, $\A(S^2_1)$ is the algebra of polynomial functions on the unit sphere and we can write the generators in terms of the standard embedding coordinates as 
\[
\begin{split}
A &= \tfrac12 (1+ z)\\
B &= \tfrac12(x+iy) \mbox.
\end{split}
\]
The commutator in $\A(S^2_q)$ gives the Poisson bracket in the commutative algebra $\A(S^2_1)$,
\[
\begin{split}
[A,B]_- &= (q^2-1)BA \quad\implies \quad
 \{A,B\} = 2iAB \\
[B,B^*]_- &= (q^2-q^{-2})A^2 - (1-q^{-2})A \quad\implies \quad
 \{B,B^*\} = 2iA(2A-1) \mbox.
\end{split}
\]
The complex coordinate $\zeta$ is related to these generators by $\zeta = A^{-1}B$, and so we can compute the Poisson bracket
\[
\{\zeta,\bar\zeta\} = 2i(\zeta\bar\zeta +1)
\mbox.\]
This agrees with \eqref{zeta.Poisson} if $c=-2$. So, we see that the Podle\'{s} standard sphere does satisfy the condition for the existence of a deformed spectral triple. It seems likely that the \Dabrowski-Sitarz example is indeed a deformed spectral triple, but this has not yet been checked explicitly.

\subsection{The Dual of a Lie Algebra}
Let $\g^*$ be the linear dual of a semisimple Lie algebra with a constant, $\g$-invariant Riemannian metric (e.~g., the Cartan-Killing form). The natural Poisson structure is given explicitly by $\pi^{ij}=C^{\;ij}_kx^k$, where $C^{\;ij}_k$ are the structure coefficients of $\g$ and $x^k$ are the coordinates. The Levi-Civita connection is simply the trivial connection given by partial derivatives on the vector space $\g^*$.

The covariant derivative of the Poisson field is simply $\nabla_k\pi^{ij} = C^{\;ij}_k$. So we have,
\[
\nabla_j \pi^{ij} = C^{\;ij}_j = 0
\mbox.\]
That is, $\g^*$ satisfies the conclusion of Thm.~\ref{Trace}. The possibility of smoothly deforming integration into a trace is not ruled out. This is discussed more extensively in \cite{b-b-g-w}.

In order to compute the curvature of the metric contravariant connection, it is sufficient to work with constant $1$-forms. The space of constant $1$-forms is naturally identified with $\g$. Let $\alpha,\beta,\gamma \in \g$. 
Let $x\cdot \alpha = x^i\alpha_i$ be the linear function on $\g^*$ with derivative $\alpha$. The definition of the Poisson structure on $\g^*$ can be stated succinctly as $\{x\cdot\alpha,x\cdot \beta\} = x\cdot [\alpha,\beta]$. With this, we see that the Koszul bracket of constant $1$-forms is simply the $\g$-bracket,
\[
[\alpha,\beta]_\pi = [d(x\cdot\alpha),d(x\cdot\beta)]_\pi = d\{x\cdot\alpha,x\cdot\beta\} = [\alpha,\beta]
\mbox.\]
We can now construct the metric contravariant connection. Because the metric inner products of $\alpha$, $\beta$ and $\gamma$ are constant over $\g^*$, eq.~\eqref{metric.connection} simplifies to 
\[
\langle\D_\alpha\beta,\gamma\rangle
= \tfrac12\left[\langle[\gamma,\alpha],\beta\rangle - \langle[\beta,\gamma],\alpha\rangle + \langle[\alpha,\beta],\gamma\rangle\right]
= \tfrac12 \langle[\alpha,\beta],\gamma\rangle
\mbox,
\]
using that $\langle[\alpha,\beta],\gamma\rangle = \langle[\beta,\gamma],\alpha\rangle = \langle[\gamma,\alpha],\beta\rangle$. So, the connection is given by
\[
\D_\alpha\beta = \tfrac12[\alpha,\beta]
\mbox.
\]
The definition \eqref{curvature} of curvature now gives,
\[
\begin{split}
K(\alpha,\beta)\gamma &= \D_\alpha\D_\beta\gamma - \D_\beta\D_\alpha\gamma - \D_{[\alpha,\beta]_\pi}\gamma \\
&= \tfrac14[\alpha,[\beta,\gamma]] - \tfrac14[\beta,[\alpha,\gamma]] - \tfrac12[[\alpha,\beta].\gamma] \\
&= -\tfrac14[[\alpha,\beta],\gamma]
\mbox.
\end{split}
\]
We see that $K\neq 0$, and so $\g^*$ does not admit a deformed spectral triple.

\subsection{Quantum Groups}
Suppose that $G$ is a semisimple Lie group with Lie algebra $\g$ and a left and right invariant Riemannian metric. Consider a deformation of $G$ into quantum groups. Because of semisimplicity, the Poisson field for such a deformation is of the form,
\beq
\label{PB.G}
\pi = r_{\mathrm{R}} - r_{\mathrm{L}}
\mbox,
\eeq
where $r_{\mathrm R}$ and $r_{\mathrm L}$ are respectively the right and left invariant bivector fields which are both equal to some $r\in \Wedge^2\g$ at the identity $e\in G$. 

The Levi-Civita connection is left and right invariant. The divergence $\nabla_j\pi^{ij}$ is a vector field on $G$ which decomposes as a sum of left and right invariant vector fields. If we identify $\g$ with left invariant vector fields on $G$, then the $\g$-bracket is identified with the Lie bracket. For $X,Y\in \g$, the Levi-Civita connection is given by
\[
\nabla_XY = \tfrac12[X,Y]
\mbox.\]
Using index notation and the structure coefficients $C^i_{\:jk}$, the covariant derivative of $r_{\mathrm L}$ is
\[
\nabla_k r^{ij}_{\mathrm L} = \tfrac12 r^{aj}C^i_{\:ka} + \tfrac12 r^{ia} C^j_{\:ka}
\mbox,
\]
and the divergence is
\beq
\label{div.r}
\nabla_j r^{ij}_{\mathrm{L}} = \tfrac12 r^{kj} C^i_{\:jk} + \tfrac12 r^{ik} C^j_{\:jk}
= -\tfrac12 r^{jk} C^i_{\:jk}
\mbox.
\eeq
This amounts to applying the Lie bracket to $\Wedge^2\g$; for example, if $r=X\wedge Y$, then \eqref{div.r} would be $-[X,Y]$.
This is the left-invariant part of $\nabla_j\pi^{ij}$. The right invariant part is essentially the same. The condition that $\nabla_j\pi^{ij}=0$ is only satisfied if \eqref{div.r} vanishes. So, integration over $G$ can only be smoothly deformed to a trace on quantum groups  if \eqref{div.r} vanishes.  For the Drinfel'd-Jimbo $r$-matrix (see, e.~g., \cite{fuc}), \eqref{div.r} is a nonzero element of the Cartan subalgebra and thus the condition is not satisfied.

The Koszul bracket of left-invariant $1$-forms is left-invariant (see \cite{wei2,l-w1,vai2}). If we identify left-invariant $1$-forms with the dual, $\g^*$, of the Lie algebra, then this bracket is given by,
\[
[\alpha,\beta]_\pi = \ad^*_{r\alpha}\beta - \ad^*_{r\beta}\alpha
\mbox,
\]
where $\ad^*$ is the coadjoint representation of $\g$ on $\g^*$ and $(r\alpha)^j := r^{ij}\alpha_i$. 

If $\alpha,\beta,\gamma\in\Omega^1(G)$ are left-invariant, the metric inner products $\langle\alpha,\beta\rangle$ \emph{et cetera} are constant functions on $G$. This simplifies eq.~\eqref{metric.connection} for the metric contravariant connection to
\begin{multline*}
\langle \D_\alpha\beta,\gamma\rangle 
= \tfrac12 \left(\langle\ad^*_{r\gamma}\alpha,\beta\rangle + \langle\ad^*_{r\gamma}\beta,\alpha\rangle  + \langle\ad^*_{r\alpha}\beta,\gamma\rangle\right. \\
\left.-\langle\ad^*_{r\alpha}\gamma,\beta\rangle - \langle\ad^*_{r\beta}\gamma,\alpha\rangle - \langle\ad^*_{r\beta}\alpha,\gamma\rangle \right)
\mbox.
\end{multline*}
Invariance of the metric shows that $\langle\ad^*_{r\alpha}\beta,\gamma\rangle + \langle\ad^*_{r\alpha}\gamma,\beta\rangle = 0$, \emph{et cetera}.
So, the metric contravariant connection on left-invariant $1$-forms is simply,
\beq
\D_\alpha\beta = \ad^*_{r\alpha}\beta
\mbox.
\eeq
Note that the metric has factored out of this formula. The connection is determined by the invariance of the metric. 

We can now compute the curvature,
\[
K(\alpha,\beta)\gamma = \ad^*_{[r\alpha,r\beta]-r[\alpha,\beta]_\pi}\gamma
= \ad^*_{[r,r](\alpha,\beta)}\gamma
\mbox,
\]
where $[r,r]\in\Wedge^3\g$ is treated as a map $[r,r] : \g^*\otimes \g^* \to \g$. This only vanishes if $[r,r]$ is constructed from the center of $\g$, but $G$ is semisimple, so the curvature only vanishes if
\beq
\label{Y-B}
[r,r] = 0
\mbox.
\eeq
This means that if a quantum group deformation of $G$ admits a deformed spectral triple, then the corresponding classical $r$-matrix must satisfy  \eqref{Y-B}, which is the well-known ``classical'' Yang-Baxter equation.

\section{Conclusions}
\label{Conclusions}
What are the implications of these results? Of course this depends on what, if anything, one wants to do with noncommutative geometry.

Suppose that space-time noncommutativity exists and is relevant to the classical physics of general relativity. We know by observation that noncommutativity can be ignored to a very good approximation, therefore we ought to be able to treat noncommutativity as a perturbation of commutative space-time. One should therefore consider general relativity with first order noncommutative corrections described by a Poisson field. However, the result is not general relativity to leading order! As I explained in Section \ref{Solutions}, consistency with a noncommutative deformation restricts the geometry of space-time. This makes space-time rigid, probably too rigid for the gravitational field to propagate in the way it is observed to. In $4$ dimensions, compatibility with a symplectic deformation leaves the metric with at most $2$ degrees of freedom per point of space-time, rather than the usual $6$ degrees of freedom. This sort of noncommutativity in the space-time of classical physics thus appears to be ruled out. 

This leaves the possibility of noncommutative space-time at some quantum level. In such a domain, there is less reason to assume that space-time is approximately commutative at all, but one can still consider the possibility. If a noncommutative deformation describes the geometry of space-time in some model ``before quantization'' or the geometry of space-time in some exotic limit of physics, then the Poisson field can be taken as a physical field and my obstructions can be taken as physical conditions (perhaps equations of motion) on  this field and the metric. 

One might also imagine that the rigidity engendered by noncommutativity is the mechanism responsible for freezing out physically unobserved extra dimensions, such as are found in Kaluza-Klein, supergravity, and string theories. A similar idea has been suggested by Dubois-Violette, Kerner, Madore, Doplicher, Fredenhagen, and Roberts \cite{d-f-r,d-k-m2}; with the motivation of preserving Poincar\'{e} symmetry, they propose a noncommutative space-time which is a deformation of a higher-dimensional manifold.

Another possibility is that my hypotheses in Theorems \ref{Trace}, \ref{Flat}, and \ref{Metric} are too restrictive. There are several ways in which this might be so.

Perhaps the assumptions of smoothness in the deformation of geometry are too restrictive. This would be surprising given how well smoothness works in the deformation of the algebra alone.

Perhaps deformations should be considered in which noncommutativity is only ``nonperturbative''. That is, it does not appear at finite order in $\k$. This is not an altogether unreasonable possibility. Consider a lattice approximation. That is, let $\M$ be a compact Riemannian manifold and $X_j\subset \M$ a sequence of finite subsets which grow uniformly dense at every point as $j\to\infty$. Construct the subset
\[
X := \left(\{0\}\times\M\right)\cup\bigcup_{j=1}^\infty \left(\{j^{-1}\}\times X_j\right) 
\; \subset [0,1]\times\M
\mbox.
\]
 Now let $\AA\subset\C(X)$ be the algebra of functions on $X$ which are restrictions of smooth functions on $[0,1]\times\M$ and let $\k\in\AA$ be the coordinate on the interval $[0,1]$. This is a smooth \emph{commutative} deformation of $\M$, but the algebra is actually unchanged to all orders in $\k$. The algebra $\AA/\k^m\AA$ is simply the tensor product of $\C^\infty(\M)$ with the algebra $\co[\k]/\k^m$. As this can happen with a commutative deformation, it can certainly happen with a noncommutative deformation. Of course, in this case there may not be anything as convenient as a Poisson field to describe the deformation.

Perhaps the structures I chose to deform are too restrictive. If one accepts the idea of a smooth deformation characterized by a Poisson field, then my results here can be taken as a guide to the structures one should not in general expect to be able to deform. If $\pi$ has nonzero divergence, then don't expect to deform integration by traces. If there does not exist a flat, torsion-free contravariant connection, then don't expect to deform $1$-forms with a gradient operator that is a derivation. If the metric contravariant connection is not flat, then don't expect to construct a well-behaved family of spectral triples.

Another possibility is the sort of noncommutativity that occurs in the almost commutative spaces of Connes-Lott models. In this case, the noncommutative space is the product of a commutative manifold and a small noncommutative space. In such a scenario, the noncommutativity is not given by a deformation. Such a space appears to be commutative not because of a scaling limit, but because the noncommutativity isn't interpreted as geometry, but as gauge and Higgs fields.

\subsection*{Acknowledgements}
I wish to thank Ludwik \Dabrowski, Lee Smolin, Ted Jacobson, and Alan Weinstein for their comments on earlier versions of this paper. This work was initiated while visiting the IH\'ES.

\providecommand{\href}[2]{#2}

\end{document}